\documentclass[12pt,leqno]{amsart}  
\usepackage{amsmath,amstext,amsthm,amssymb,amsxtra}
\usepackage[colorlinks,citecolor=red,pagebackref,hypertexnames=false]{hyperref}
\usepackage[backrefs]{amsrefs}
 
\usepackage{pgf,pgfarrows}

\usepackage{amsmath, amssymb, amscd, amssymb}
\usepackage{latexsym}

\begin{document}

\newtheorem{theorem}[equation]{Theorem}    
\newtheorem{proposition}[equation]{Proposition}
\newtheorem{conjecture}[equation]{Conjecture}
\def\theconjecture{\unskip}
\newtheorem{corollary}[equation]{Corollary}
\newtheorem{lemma}[equation]{Lemma}
\newtheorem{observation}[equation]{Observation}
\theoremstyle{definition}
\newtheorem{definition}[equation]{Definition}
\newtheorem{remark}[equation]{Remark}
\newtheorem*{Acknowledgment}{Acknowledgment}

\numberwithin{equation}{section}

\def\abs#1{\lvert#1\rvert}
\def\Abs#1{\bigl\lvert#1\bigr\rvert}
\def\ABs#1{\Bigl\lvert#1\Bigr\rvert}
\def\ABS#1{\biggl\lvert#1\biggr\rvert}
\def\ABSS#1{\Biggl\lvert#1\Biggr\rvert}

\def\name #1 #2{\operatorname {\mathsf #1}(#2)}

\newcommand{\ann}{{\rm scl}}
\newcommand{\scl}{{\rm scl}}
\newcommand{\Sc}{{\mathcal {S}}}
\newcommand{\Sb}{{\mathbf{S}}}

\def\norm#1.#2.{\lVert#1\rVert_{#2}}
\def\Norm#1.#2.{\bigl\lVert#1\bigr\rVert_{#2}}
\def\NOrm#1.#2.{\Bigl\lVert#1\Bigr\rVert_{#2}}
\def\NORm#1.#2.{\biggl\lVert#1\biggr\rVert_{#2}}
\def\NORM#1.#2.{\Biggl\lVert#1\Biggr\rVert_{#2}}


\def\ip#1,#2,{\langle #1,#2\rangle}
\def\Ip#1,#2,{\bigl\langle#1,#2\bigr\rangle}
\def\IP#1,#2,{\Bigl\langle#1,#2\Bigr\rangle}

\def\abs#1{\lvert#1\rvert}
\def\Abs#1{\bigl\lvert#1\bigr\rvert}
\def\ABs#1{\Bigl\lvert#1\Bigr\rvert}
\def\ABS#1{\biggl\lvert#1\biggr\rvert}
\def\ABSS#1{\Biggl\lvert#1\Biggr\rvert}

\def\eqdef{\stackrel{\mathrm{def}}{=}}


\title[Lipschitz Kakeya Maximal Function]{On a Lipschitz  Variant  of the  Kakeya Maximal Function  }

\author[M. Lacey \& Xiaochun Li]{Michael Lacey and Xiaochun Li}

\address{Michael Lacey \\ School of Mathematics \\ 
Georgia Institute of Technology \\ Atlanta GA 30332 }

\address{        Xiaochun Li\\
         	Department of Mathematics\\
        University of Illinois\\
        Urbana IL 61801}

\email{xcli@math.uiuc.edu}
\thanks{The authors are supported in part by NSF grants.  M.L. is a Guggenheim Fellow.}

%
%

%

\maketitle

\section{The Lipschitz Kakeya Maximal Functions}

 In this paper, 
we are concerned with a maximal function estimate over rectangles in the plane, which 
have a prescribed maximal length, but arbitrary orientation and width.  That is, we are concerned 
with a certain variant of the Kakeya maximal function. 
The rectangles used in the maximal function will be specified in a particular way by 
a  Lipschitz choice of directions in the plane. 

Our results are  stated first; afterwords we place it in the context of the literature. 
We set notations and conventions.  A \emph{rectangle} is determined as follows.  Fix a choice 
of unit vectors in the plane $ (\operatorname e,\operatorname e ^{\perp})$, with $ \operatorname 
e ^{\perp}$ being the vector $ \operatorname e$ rotated by $ \pi /2$.  Using these vectors 
as coordinate axes, a rectangle is a product of two intervals $ R= I\times J$.  We will 
insist that $ \abs{ I}\ge \abs{ J}$, and use the notations 
\begin{equation}\label{e.LW}
\name L R=\abs{ I},\qquad \name W R=\abs{ J}
\end{equation}
for the length and width respectively of $ R$. 
 
The \emph{interval of uncertainty of  $ R$} is the subarc  $ \name EX R $ of the unit circle in the plane, 
centered at $ \operatorname e$, and of length $ \name W R/\name L R$.  See Figure~1.

 \begin{figure}
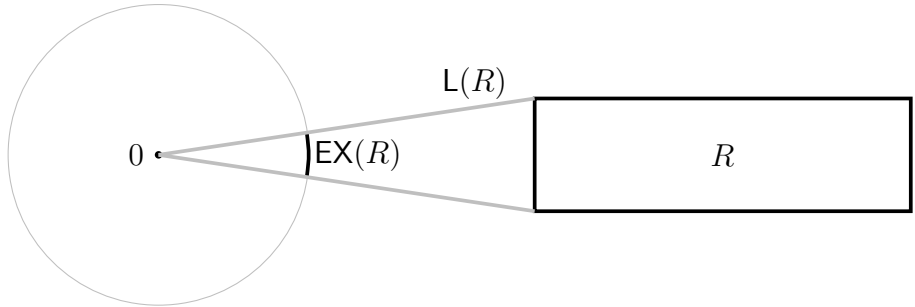

	 \begin{pgfpicture}{0cm}{0cm}{8cm}{4.5cm}
	 \begin{pgftranslate}{\pgfxy(-.15,2)}
	{\color{lightgray}  
		\pgfcircle[stroke]{\pgfxy(0,0)}{2cm}
	}
	\pgfcircle[fill]{\pgfxy(0,0)}{0.05cm}
	\pgfputat{\pgfxy(-.3,0)}{\pgfbox[center,center]{$0$}}
	\pgfputat{\pgfxy(2.65,0)}{\pgfbox[center,center]{$\mathsf{ EX}(R) $}}
	\pgfputat{\pgfxy(7.5,0)}{\pgfbox[center,center]{$R $}}
	\pgfsetlinewidth{1.4pt}
	\pgfrect[stroke]{\pgfxy(5,-.75)}{\pgfxy(5,1.5)}
	\pgfmoveto{\pgfxy(1.97,.3)}
	\pgfcurveto{\pgfxy(2.005,0)}{\pgfxy(2.005,0)}{\pgfxy(1.97,-.3)}\pgfstroke
			\pgfputat{\pgfxy(4.2,.95)}{\pgfbox[center,center]{$ \mathsf L (R)$}}
{\color{lightgray} 
			\pgfmoveto{\pgfxy(0,0)} 
			\pgflineto{\pgfxy(5,.75)} 	\pgfstroke
			\pgfmoveto{\pgfxy(0,0)} 
			\pgflineto{\pgfxy(5,-.75)}	\pgfstroke
	}
	\end{pgftranslate}
	\end{pgfpicture}
	\caption{ An example eccentricity interval $ \mathsf {EX}(R)$.  The circle on 
	the left has radius one.}
	\end{figure}

We now fix a Lipschitz map $ v$ of the plane into the unit circle. 
We only consider rectangles $ R$ with 
\begin{equation} \label{e.shortlength}
\name L R\le{} (100 \norm v. \text{Lip} .) ^{-1}\,.
\end{equation}  
For such a rectangle $ R$, set $ \name V R=R\cap v ^{-1}(\name EX R )$. 

For $ 0<\delta <1$, we consider the maximal functions 
\begin{equation}\label{e.Mdef}
\operatorname M _{v,\delta }f(x) \eqdef \sup _{\substack{ \abs{ \name V R }\ge \delta \abs{ R} }}
\frac {\mathbf 1 _{R}(x)} {\abs{ R} } \int _{R} \abs{f(y)}\; dy .
\end{equation}
That is we only form the supremum over rectangles for which the vector field
 lies in the interval of uncertainty for a fixed positive proportion $ \delta $ of 
 the rectangle.

\begin{theorem}\label{t.main}
The maximal function $\operatorname M_{\delta,v}$ is bounded from $L^2({\mathbb R}^2)$
to $L^{2, \infty}({\mathbb R}^2)$ with norm at most $ \lesssim \delta ^{-1/2}$. That is, 
for any $\lambda >0$,  and $ f\in L^2 (\mathbb R^2)$, this inequality holds: 
\begin{equation}\label{max1}
 \lambda ^{-2}
 \abs{\{x\in {\mathbb R}^2: \operatorname M_{\delta,v}f(x)>\lambda \}} 
 \lesssim \delta^{-1}\lambda^{-2}\norm f.2.^2\,.
\end{equation}
The norm estimate in particular is independent of the Lipschitz vector field $ v$. 
\end{theorem}

The Lipschitz assumption is sharp; it cannot be 
weakened to H\"older continuity of any index less than one. 

\bigskip

The Kakeya maximal function is typically defined as 
\begin{equation*}
\operatorname M _{\textup K, \epsilon } f(x)=\sup _{\abs{ \name EX R }\ge \epsilon } 
\frac {\mathbf 1 _{R}(x) } {\abs{ R}} \int _{R} \abs{f(y)}\; dy \,, \qquad \epsilon >0.
\end{equation*}
One is forced to take $ \epsilon >0$ due to the existence of the Besicovitch set.  It is 
a critical fact that the norm of this operator admits norm bound on $ L^2$ that is 
logarithmic in $ \epsilon ^{-1}$.  See  C{\'o}rdoba and Fefferman \cite{MR0476977}, 
and Str{\"o}mberg \cites {MR0487260
,
MR0481883}.  Subsequently, there have been several refinements of this observation, 
we cite only Nets H. Katz \cite {nets}, 
Alfonseca, Soria and Vargas \cite{MR1960122}, 
and Alfonseca \cite{MR1942421}.  These papers contain additional references. 

Our variant is inspired by questions of Zygmund and E.M.~Stein concerning certain 
degenerate Radon transforms.  For Lipschitz vector fields $ v$, consider the maximal 
function and Hilbert transforms along line segments determined by $ v$.  For 
$ \nu =(100 \norm v.\textup{Lip}.) ^{-1}$, set 
\begin{align*}	
\operatorname M_v f(x) &\eqdef \sup _{0<t<\nu } (2t) ^{-1} 
\int _{ -t } ^{t } \abs{ f(x-yv(x)) }\; dy\,, 
\\   
\operatorname H _v f(x) & \eqdef \text{p.v.}\int _{-\nu } ^{\nu } f(x-y v (x))\frac {dy}y\,.
\end{align*}
The question of the boundedness of the $\operatorname  M_v$ on e.g.~$ L^2(\mathbb R^2)$ is attributed 
to Zygmund, motivated by   constructions of the Besicovitch set which
would show that the Lipschitz condition 
is sharp.  E.M.~Stein raised a similar question about the Hilbert transform analog \cite
{stein}. 

Both operators are examples of Radon transforms.  And that theory generally applies under 
suitable additional geometric conditions placed on the vector field $ v$.  Thus, 
there are positive results in the case that the vector field is analytic, due to 
Nagel, Stein and Wainger \cite{MR81a:42027}.  Bourgain \cite{bourgain} provided an extension to the case of 
real analytic vector fields, which exhibit a range of degeneracies that analytic 
vector fields can't possess.  The theory of Radon transforms has a beautiful exposition in 
the article of Christ, Nagel, Stein and Wainger \cite{MR2000j:42023}.  Another paper 
relevant to the Radon transforms of this paper 
is one by Carbery, Seeger and Wainger, Wright \cite{carbery}.  Both papers contain a 
wide set of references to the literature on Radon transforms.

Nets Katz \cite{MR1979942} provided a partial result in the direction of Zygmund's conjecture.

We speculate that the norm bounds obtained in the Theorem  is optimal. 
Namely, that the operator $ \operatorname M _{v,\delta }$ will in general 
have a weak type bound on $ L^2$ that is at least as big as $c \delta ^{1/2}$, and is 
unbounded on $ L^p$ for $ 1<p<2$.

\bigskip 

The notation $ A \lesssim B$ means that $ A\le KB$ for an absolute, but unspecified, 
constant $ K$.  $ A\simeq B$ means that $ A \lesssim B \lesssim A$.  $ \mathbf 1 _{A}$ 
means the indicator function of the set $ A$.

\section{Proof of Theorem~\ref{t.main}} 

\subsection*{The Covering Lemma Conditions} 

We adopt the covering lemma approach of 
C{\'o}rdoba and R.~Fefferman \cite{MR0476977}.  To this end, we regard the choice of vector field 
$ v$ and $ 0<\delta <1$ as fixed.  Let $ \mathcal R$ be any finite collection of rectangles 
obeying the conditions (\ref{e.shortlength}) and $ \abs{ \name V R }\ge \delta \abs{ R}$. 
We show that $ \mathcal R$ has a decomposition into disjoint collections $ \mathcal R'$ 
and $ \mathcal R''$ for which these estimates hold.   
\begin{align}\label{e.2<1}
\NOrm \sum _{R\in \mathcal R'} \mathbf 1 _{R}.2.^2 &\lesssim \delta ^{-1 } 
\NOrm \sum _{R\in \mathcal R'} \mathbf 1 _{R}.1.\,,
\\  \label{e.bigcup}
\ABs{\bigcup _{R\in \mathcal R''}R} &\lesssim 
\NOrm \sum _{R\in \mathcal R'} \mathbf 1 _{ R}.1.
\end{align}
The first of these conditions is the stronger one, as it bounds the $ L^2$ norm squared 
by the $ L^1$ norm; the verification of it will occupy most of the proof. 

\smallskip 

Let us see how to deduce Theorem~\ref{t.main}.  Take $ \lambda >0$ and $ f\in L^2$  which is non negative and 
 of norm one. Set $ \mathcal R$ to be all the rectangles  $ R$ 
 of prescribed maximum length as given in (\ref{e.shortlength}), 
 density with respect to the vector field, namely $ \abs{ \name V R}\ge \delta \abs{ R}$,  and  
\begin{equation*}
\int _{R} f(y)\; dy\ge{} \lambda \abs{ R} \,.
\end{equation*}
We should verify the weak type inequality 
\begin{equation} \label{e.wweak}
\lambda \ABs{\bigcup _{R\in \mathcal R}R} ^{1/2} \lesssim \delta ^{-1/2}\,.
\end{equation}

Apply the decomposition to $ \mathcal R$.  Observe that 
\begin{align*}
\lambda \NOrm \sum _{R\in \mathcal R'} \mathbf 1 _{R}.1. 
&\le \IP f, \sum _{R\in \mathcal R'} \mathbf 1 _{R},
\\
&\le \NOrm \sum _{R\in \mathcal R'} \mathbf 1 _{R}.2.
\\
& \lesssim \delta ^{-1/2} \NOrm \sum _{R\in \mathcal R'} \mathbf 1 _{R}.1. ^{1/2}\,.
\end{align*}
Here of course we have used (\ref{e.2<1}).  This implies that 
\begin{equation*}
\lambda \NOrm \sum _{R\in \mathcal R'} \mathbf 1 _{R}.1. ^{1/2} \lesssim  \delta ^{-1/2}.
\end{equation*}
Therefore clearly (\ref{e.wweak}) holds for the collection $ \mathcal R'$. 

Concerning the collection $ \mathcal R''$, apply (\ref{e.bigcup}) to see that  
\begin{align*}
\lambda \ABs{\bigcup _{R\in \mathcal R''}R} ^{1/2}  
\lesssim \lambda \NOrm \sum _{R\in \mathcal R'} \mathbf 1 _{R}.1. ^{1/2}
\lesssim \delta ^{-1/2}\,.
\end{align*}
This completes our proof of (\ref{e.wweak}).

\smallskip  

The remainder of the proof is devoted to  the proof of (\ref{e.2<1}) and (\ref{e.bigcup}).

\subsection*{The Covering Lemma Estimates} 

\subsubsection*{Construction of $ \mathcal R'$ and $ \mathcal R''$}
In the course of the proof, we will need several recursive procedures.  The 
first of these occurs in the selection of $ \mathcal R'$ and $ \mathcal R''$. 

We will have need of one large constant $ \kappa $, of the order of say $ 100$, but whose 
exact value does not concern us. Using this notation hides 
distracting terms.

Let $\operatorname M_{\kappa}$ be a maximal function given as 
\begin{equation*}
\operatorname M _{\kappa } f(x)=\sup _{s>0} \max \Bigl\{ 
 s ^{-2}\int _{x+sQ} \abs{ f(y)} \; dy\, ,\ \sup _{\omega \in \Omega } s ^{-1 }
 \int _{-s} ^{s} \abs{ f(x+ \sigma \omega  )} \; d \sigma \Bigr\}\,.
\end{equation*}
Here, $ Q$ is the unit cube in plane, and $ \Omega $ is a set of uniformly distributed 
points on the unit circle of cardinality equal to $ \kappa $. 
It follows from the usual weak type bounds 
that this operator maps $ L^1(\mathbb R^2)$ into weak $ L^1(\mathbb R^2)$. 

\smallskip 

To initialize the recursive procedure, set 
\begin{align*}
\mathcal R'&\leftarrow \emptyset\,,
\\
\mathsf{STOCK} & \leftarrow \mathcal R\,.
\end{align*}

The main step is this while loop.  While $ \mathsf{STOCK}$ is not empty, 
select $ R \in \mathsf{STOCK}$ subject to the criteria that first it have a 
maximal  length $ \mathsf L (R)$, and second that it have minimal value 
of $ \abs{ \name EX R}$. 
Update 
\begin{equation*}
\mathcal R'\leftarrow \mathcal R'\cup\{R\}. 
\end{equation*}
Remove $ R$ from $ \mathsf{STOCK}$. As well, remove any rectangle $ R'\in \mathsf{STOCK}$
which is also contained in 
\begin{equation*}
\Bigl\{ \operatorname M_{\kappa} \sum _{R\in \mathcal R'} \mathbf 1 _{\kappa R}\ge{} \kappa ^{-1}
\Bigr\}.
\end{equation*}

As the collection $ \mathcal R$ is finite, the while loop will terminate, and at this point 
we set $ \mathcal R''{} \eqdef  \mathcal R- \mathcal R'$. 
In the course of the argument below, we will refer the order in which rectangles 
where added to $ \mathcal R'$. 

\smallskip

With this construction, it is obvious that (\ref{e.wweak}) holds, with a bound 
that is a function of $ \kappa $.  Yet, $ \kappa $ is an absolute constant, so 
this dependence does not concern us. And so the 
rest of the proof is devoted to the verification of (\ref{e.2<1}). 

\smallskip

An important aspect of the qualitative nature of the interval of eccentricity is 
encoded into this algorithm.  We will choose $ \kappa $ so large that this is true: 
Consider two rectangles $ R$ and $ R'$ with $ R\cap R'\neq\emptyset$, 
$ \name L R\ge \name L {R'}$, $ \name W R\ge \name W {R'} $,  
$ \abs{\name {EX} {R}}\le{} \abs{\name {EX} {R'} }$ and $\name EX R\subset 10\name EX {R'} $
then we have    
\begin{equation}\label{e.ex==>}
R'\subset \kappa R\,.
\end{equation}
See Figure~2.

 \begin{figure}  
\begin{center}
  \begin{pgfpicture}{0cm}{0cm}{8cm}{4cm}
  \begin{pgftranslate}{\pgfxy(4,2.5)}
	\pgfrect[stroke]{\pgfxy(-3,-.3)}{\pgfxy(6,.6)}
	\pgfputat{\pgfxy(3.25,0)}{\pgfbox[center,center]{$ R$}}
	\begin{pgfrotateby}{\pgfdegree{20}}
	\pgfrect[stroke]{\pgfxy(-2,-.07)}{\pgfxy(4,.14)}
	\end{pgfrotateby}
	\pgfputat{\pgfxy(1.8,1.1)}{\pgfbox[center,center]{$ R'$}}
 {\color{lightgray} 
 	\pgfrect[stroke]{\pgfxy(-5.5,-.8)}{\pgfxy(11,1.6)}
	\pgfputat{\pgfxy(4.7,1.1)}{\pgfbox[center,center]{$\kappa R$}}
}
	\end{pgftranslate}
	\end{pgfpicture}
	\end{center} 
	\caption{} 
	\end{figure}

\subsection*{Uniform Estimates}

We estimate the left hand side of (\ref{e.2<1}).  In so doing we expand the square, 
and seek certain uniform estimates.  Expanding the square on the left hand side 
of (\ref{e.2<1}), we can estimate 
\begin{equation*}
\textup{l.h.s.\thinspace of } (\ref{e.2<1})
\le \sum _{R\in \mathcal R'} \abs{ R}+2 \sum _{(\rho ,R)\in \mathcal P} \abs{ \rho \cap R}
\end{equation*}
where $ \mathcal P$ consists of all pairs $ (\rho ,R)\in \mathcal R' \times \mathcal R'$ 
such that $ \rho \cap R\neq \emptyset$, and 
$ \rho $ was selected to be a member of $ \mathcal R'$ before 
$ R$ was.  It is then automatic that $ \mathsf L (R )\le \mathsf L (\rho )$.
And since the density of all tiles is positive, it follows that 
$ \operatorname {dist} (\mathsf {EX} (\rho ), \mathsf {EX} (R))
\le 2 \norm v.\textup{Lip}. \mathsf L (\rho )<\tfrac 1 {50}$.

We will split up the collection $ \mathcal P$ into 
subcollections $ \{\mathcal S _{R} \mid R\in \mathcal R'\}$ 
and $ \{\mathcal T _{\rho }\mid \rho \in \mathcal R'\}$.

\medskip 

For a rectangle $ R\in \mathcal R'$, we take 
$ \mathcal S _{R}$ to consist of all rectangles $ \rho $ such that 
(a) $ (\rho ,R)\in \mathcal P$;  and (b) $\name EX \rho \subset 
10\name EX R $. 
We assert that 
\begin{equation}\label{e.uni1}
\sum_{\rho \in \mathcal S_R} \abs{ R\cap \rho } \le 
\abs{ R},\qquad R\in \mathcal R. 
\end{equation}

This estimate is in fact easily available to us. 
Since the rectangles $ \rho \in \mathcal S_R$ were selected to be 
in $ \mathcal R'$ before $ R$ was, we cannot have the inclusion 
\begin{equation} \label{e.Xuni}
R\subset \Bigl\{ \operatorname M_{\kappa} \sum _{\rho \in \mathcal S_R} \mathbf 1 _{\kappa \rho } 
> \kappa ^{-1} \Bigr\}\,.
\end{equation}
Now the  rectangle $ \rho $ are  also longer.  Thus, if (\ref{e.uni1}) 
 does not hold,
 we would compute the maximal function of 
 \begin{equation*}
\sum _{\rho \in \mathcal S_R} \mathbf 1 _{\kappa \rho }
\end{equation*}
in a direction which is close, within an error of $ 2 \pi /\kappa $, of being orthogonal 
to the long direction of $ R$.  In this way, we will contradict (\ref{e.Xuni}).

\medskip 

The second uniform estimate that we need is as follows.  
For fixed $ \rho $, set $ \mathcal T _{\rho }$ to be the set of all rectangles 
$ R$ such that 
(a) $ (\rho ,R)\in \mathcal P$ and (b) 
$\name EX \rho \not\subset  10\name EX R $. 
We assert that 
\begin{equation}\label{e.uni2}
\sum_{R\in \mathcal T _{\rho }} 
\abs{ R\cap \rho } \lesssim \delta ^{-1} \abs{ \rho },\qquad \rho \in \mathcal R'. 
\end{equation}
This proof of this inequality is more involved, and taken up in the next subsection. 

\begin{remark}\label{r.rho} In the proof of (\ref{e.uni2}), it is not 
necessary that $ \rho \in \mathcal R'$.  Writing $ \rho =I _{\rho } \times J _{\rho }$, 
in the coordinate basis $ \operatorname e$ and $ \operatorname e _{\perp}$, 
we could take any rectangle of the form $ I \times J _{\rho }$. 
\end{remark}

\medskip 

These two estimates conclude the proof of (\ref{e.2<1}). 
For any two distinct rectangles $\rho, R\in \mathcal P$,
we will have either $ \rho \in \mathcal S_R$ or $ R\in \mathcal T _{\rho }$.
Thus (\ref{e.2<1}) follows by summing (\ref{e.uni1}) on $ R$ and (\ref{e.uni2})
on $ \rho $. 

\subsection*{The Proof of (\ref{e.uni2})} 
%
%
%
%
%
%
%
%
%
%

We fix $ \rho $, and begin by making a decomposition of the collection $ \mathcal T _{\rho }$. 
Suppose that the coordinate axes for $ \rho $ are given by $ \operatorname e _{\rho }$, 
associated with the long side of $ R$, and 
$ \operatorname e ^{\perp} _{\rho }$, with the short side.  Write the rectangle as a product of intervals $ I _{\rho }\times J$, 
where $ \abs{ I _{\rho }}=\name L \rho $.  
Denote one of the endpoints of $ J$ as $ \alpha $.  
See Figure~3.

\def\thinrect#1#2#3#4{
	\pgfrect[stroke]{\pgfpoint{#3cm}{#4cm}}{\pgfpoint{#2cm}{.3cm}}
	 }
%
%
 \begin{figure}
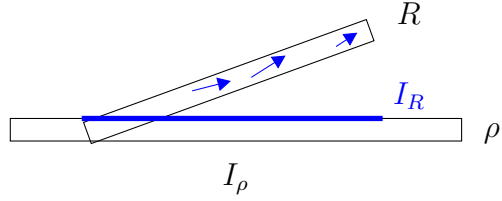
  
\begin{center} 
 \begin{pgfpicture}{0cm}{0cm}{8cm}{4cm}
\begin{pgftranslate}{\pgfxy(2,1)}
{\thinrect 0600}
 \pgfputat{\pgfxy(6.4,0.1)}{\pgfbox[center,center]{$\rho $}}
  \pgfputat{\pgfxy(3,-0.5)}{\pgfbox[center,center]{$I_\rho  $}}
 				\pgfsetendarrow{\pgfarrowtriangle{4pt}}
  		\begin{pgfrotateby}{\pgfdegree{20}}
	\thinrect 041{-.4}
	{\textcolor{blue}{
			\pgfmoveto{\pgfxy(4.5,-.3)}\pgflineto{\pgfxy(4.8,-.23)}\pgfstroke
			\pgfmoveto{\pgfxy(3.3,-.3)}\pgflineto{\pgfxy(3.8,-.19)}\pgfstroke
			\pgfmoveto{\pgfxy(2.5,-.2)}\pgflineto{\pgfxy(3,-.25)}\pgfstroke}}
	\end{pgfrotateby}
	\pgfputat{\pgfxy(5.3,1.7)}{\pgfbox[center,center]{$R$}}
	\pgfclearendarrow
	\textcolor{blue}{  \pgfsetlinewidth{2pt}
 \pgfmoveto{\pgfxy(.95,.3)}\pgflineto{\pgfxy(4.95,.3)}\pgfstroke
 \pgfputat{\pgfxy(5.3,.6)}{\pgfbox[center,center]{$I_R$}}
}
 \end{pgftranslate}  
 \end{pgfpicture}
 \end{center} 
 \caption{Notation for the proof of (\ref{e.uni2}).}
 \end{figure}

For rectangles $ R\in \mathcal T _{\rho }$, let $ I_R$ denote the orthogonal projection 
$ R$ onto the line segment $ 2I _{\rho }\times \{\alpha \}$.  Subsequently, we will consider 
different subsets of this line segment.  The first of these is as follows. 
For $ R\in \mathcal T _{\rho }$, let $ \mathsf V_R$ be the projection 
of the set $ \name V R $ onto $2I\times \{\alpha \}$.  We have 
\begin{equation} \label{e.B}
\name L R\le \abs{ I_R}\le 2 \name L R ,
\quad\text{and}\quad
\delta \name L R\lesssim  \abs{ \mathsf V_R}. 
\end{equation}

A recursive mechanism is used to decompose $ \mathcal T _{\rho }$.  Initialize 
\begin{align*}
\mathsf{STOCK} & \leftarrow \mathcal T_\rho \,,
\\
\mathcal U & \leftarrow \emptyset \,.
\end{align*}
While $ \mathsf{STOCK}\neq\emptyset$ select $ R\in \mathsf{STOCK}$ of maximal 
length. Update 
\begin{equation}\label{e.UU}
\begin{split}
\mathcal U&\leftarrow \mathcal U\cup\{R\},
\\
\mathcal U(R) &\leftarrow \{R'\in \mathsf{STOCK}\mid \mathsf V_R\cap \mathsf V_{R'}\neq\emptyset\}. 
\\
\mathsf{STOCK}&\leftarrow\mathsf{STOCK}-\mathcal U(R). 
\end{split}
\end{equation}
When this while loop stops, it is the case that $ \mathcal T_ \rho =\bigcup _{R\in \mathcal U} 
\mathcal U(R)$.

\smallskip

With this construction, the sets $ \{ \mathsf V_R \mid R\in \mathcal U\} $ are 
disjoint.  By (\ref{e.B}), we have 
\begin{equation}\label{e.U}
\sum _{R\in \mathcal U} \name L R \lesssim \delta ^{-1} \name L \rho \,.
\end{equation}
The main point, is then to verify the uniform estimate 
\begin{equation}\label{e.uni3}
\sum _{R'\in \mathcal U(R)} \abs{ R'\cap \rho } 
\lesssim  \name L R\cdot \name W \rho  \,, 
\qquad R\in \mathcal U. 
\end{equation}
Note that both estimates immediately imply (\ref{e.uni2}).

\subsubsection*{Proof of (\ref{e.uni3})}

There are three important, and more technical, facts to observe about the collections 
$ \mathcal U(R)$.

\def\thinrect#1#2#3#4{
	\pgfrect[stroke]{\pgfpoint{#3cm}{#4cm}}{\pgfpoint{#2cm}{.3cm}}
	 }

\begin{figure}
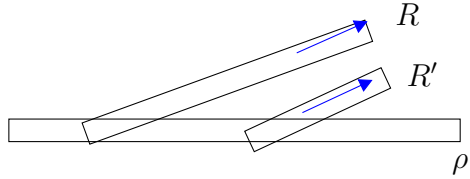
   
\begin{center} 
 \begin{pgfpicture}{0cm}{0cm}{7.5cm}{4cm}
\begin{pgftranslate}{\pgfxy(1,1)}
 \thinrect 0600
 \pgfputat{\pgfxy(6,-.3)}{\pgfbox[center,center]{$\rho $}}
 \pgfsetendarrow{\pgfarrowtriangle{4pt}}
 	\begin{pgfrotateby}{\pgfdegree{20}}
	\thinrect 041{-.4}
	\textcolor{blue}{\pgfmoveto{\pgfxy(4,-.2)}\pgflineto{\pgfxy(5,-.13)}\pgfstroke}
			\end{pgfrotateby}
	\pgfputat{\pgfxy(5.3,1.7)}{\pgfbox[center,center]{$R$}}
		 	\begin{pgfrotateby}{\pgfdegree{25}}
	\thinrect 02{2.9}{-1.5}
	\textcolor{blue}{\pgfmoveto{\pgfxy(3.7,-1.3)}\pgflineto{\pgfxy(4.7,-1.3)}\pgfstroke}
		\end{pgfrotateby}
		\pgfputat{\pgfxy(5.5,.9)}{\pgfbox[center,center]{$R'$}}
 \end{pgftranslate}  
 \end{pgfpicture}
 \end{center} \caption{Proof of Lemma~\ref{l.samedirection}:  The rectangles $R,R'\in \mathcal U(\rho ) $, 
 and so the angles $ R$ and $ R'$ form with $ \rho $ are nearly the same. }
 \end{figure}

For any rectangle $ R'\in \mathcal U(R)$, denote it's coordinate 
axes as $ \operatorname e _{R'}$ and $ \operatorname e _{R'} ^{\perp}$, 
associated to the long and short sides of $ R'$ respectively.

\begin{lemma}\label{l.samedirection}
For any rectangle $ R'\in \mathcal U (R)$ we have  
\begin{equation*}
\abs{ \operatorname e _{R'}-\operatorname e _{R}} \le \tfrac 12 
\abs{ \operatorname e _{\rho }-\operatorname e _{R}}
\end{equation*}
\end{lemma}

\begin{proof}
There are by construction, points $ x\in \name V R$ and $ x'\in  \name V {R'} $ which get projected 
to the same point on the line segment $ I_\rho \times \{\alpha \}$. 
See Figure~4. Observe that 
\begin{align*}
\abs{ \operatorname e _{R'}-\operatorname e _{R}}&\le
\abs{ \name EX {R'} }+\abs{\name  EX R}+\abs{ v(x')-v(x)}
\\
&\le 
\abs{ \name EX {R'} }+\abs{\name  EX R}+\norm v.\text{Lip}. \cdot \name L R \cdot \abs{ \operatorname e _{\rho }-\operatorname e _{R}}
\\
&\le{} \abs{ \name EX {R'} }+\abs{\name  EX R}+\tfrac1{100} \abs{ \operatorname e _{\rho }-\operatorname e _{R}}
\end{align*}
Now, $ \abs{ \name EX R}\le{} \tfrac1{5} \abs{ \operatorname e _{\rho }-\operatorname e _{R}}$, 
else we would have $ \rho \in \mathcal S _{R}$. 
Likewise, $\abs{ \name EX {R'} }\le{} \tfrac1{5} \abs{ \operatorname e _{R' }-\operatorname e _{R}} $.
And this proves the desired inequality. 
 
\end{proof}

\begin{lemma}\label{l.geo}
Suppose that there is an interval $I'\subset I _{\rho }$ such that
\begin{equation}\label{1}
\sum_{\substack{R'\in \mathcal U(R )\\ \name L {R'} \geq 8\abs I}} \abs{R'\cap I\times J} 
\geq     \abs{I\times J}\,.
\end{equation}
Then 
there is no $R''\in \mathcal U(R )$ such that $ \name L {R''} < \abs{I}$ and 
$  R''\cap 4I\times J\neq\emptyset  $. 
\end{lemma}

\begin{proof}
There is a natural angle $ \theta $ between the rectangles $ \rho $ and $ R$, 
which we can assume is positive, and is 
given by $ \abs{\operatorname e _{\rho }-\operatorname e _{R}}$.  
Notice that we have $ \theta\ge 10\abs{\name EX R} $, else we would have $ \rho \in \mathcal S _R$, 
which contradicts our construction.

Moreover, there is an important consequence of the previous Lemma~\ref{l.samedirection}: For any 
$ R'\in \mathcal U(R)$, there is a natural angle $ \theta '$ between $ R'$ and $ \rho $. 
These two angles are close.  For our purposes below, these two 
angles can be regarded as the same.

For any $ R'\in \mathcal U(R)$, we will have 
\begin{align*}
\frac{\abs{ \kappa R'\cap \rho }} {\abs{I_\rho \times J}} 
&\simeq \kappa \frac{ \name W {R'} \cdot \name W \rho }{\theta \abs{ I_\rho } \name W \rho  }
\\
&= \kappa \frac{\name W {R'} }{\theta \cdot  \abs {I_\rho } }.
\end{align*}

Recall $\operatorname M_{\kappa}$ is larger than the maximal function over $ \kappa $ uniformly 
distributed directions. Choose a direction $\operatorname e'$ from this set of  $ \kappa $ directions
that is closest to $\operatorname e ^\perp _{\rho }$. Take a line segment $\Lambda  $ in
direction $\operatorname e'$  of length  $\kappa \theta \abs{ I} $, 
and the center of $ \Lambda $ is in $4I\times J$.  See Figure~5. 
Then we have 
\begin{align*}
\frac{\abs{\kappa R'\cap  \Lambda }}{\abs{\Lambda }}& \geq 
\frac{ \name W {R'}   } {  \theta\cdot  \abs I  }
\end{align*}
Thus by our assumption (\ref{1}),  
$$
\frac1 {\abs{ \Lambda }} \sum _{R'\in \mathcal U(R)} \abs{ R'\cap \Lambda }\ge 1 \,.
$$
That is, any of the lines $ \Lambda $ are contained in the set 
\begin{equation*}
\Bigl\{  \operatorname M _{\kappa }\sum _{R\in \mathcal R'} \mathbf 1 _{R'}>\kappa ^{-1} \Bigr\}.  
\end{equation*}

\smallskip

\def\thinrect#1#2#3#4{
	\pgfrect[stroke]{\pgfpoint{#3cm}{#4cm}}{\pgfpoint{#2cm}{.3cm}}
	 }

 \begin{figure}   
\begin{center} 
\begin{pgfpicture}{0cm}{0cm}{5cm}{4cm}
\begin{pgftranslate}{\pgfxy(1.5,0)}
%
{\color{gray} \pgfrect[stroke]{\pgfxy(-3.3,0.5)}{\pgfxy(6.6,.2)} }
\pgfputat{\pgfxy(-2.7,1)}{\pgfbox[center,center]{$4I\times J$}}
\pgfrect[stroke]{\pgfxy(-1.3,0.5)}{\pgfxy(2.3,.2)}
\pgfputat{\pgfxy(.8,.2)}{\pgfbox[center,center]{$I\times J$}}
\pgfsetlinewidth{0.85pt}
\pgfline{\pgfxy(2.1,-0.7)}{\pgfxy(2.1,2.0)}
\pgfputat{\pgfxy(2.45,0)}{\pgfbox[center,center]{$\Lambda $}}
\begin{pgftranslate}{\pgfxy(-1,.5)}
\begin{pgfrotateby}{\pgfdegree{15}}
\pgfrect[stroke]{\pgfxy(-3,-.1)}{\pgfxy(8,.2)}
\end{pgfrotateby}
\end{pgftranslate}
\pgfputat{\pgfxy(3.5,1.4)}{\pgfbox[center,center]{$\kappa R'$}}
\end{pgftranslate}
\end{pgfpicture}
\end{center}
\caption{} 
\end{figure}

Clearly our construction does not permit 
any rectangle $ R''\in \mathcal U(R) $ contained in this set. 
To conclude the proof of our Lemma, we seek a contradiction.  Suppose that there 
is an $ R''\in \mathcal U(R)$ with $ \name L {R''}<\abs{ I}$ and $ R''$ intersects 
$ 2I \times J$.  The range of line segments $ \Lambda  $ we can permit is however quite 
broad.  The only possibility permitted to us is that the rectangle $ R''$ is quite 
wide.  We must have 
\begin{equation*}
\name W {R''}\ge{}\tfrac 14 \abs{  \Lambda  } =\tfrac \kappa 4 \cdot \theta \cdot  \abs I .  
\end{equation*}
This however forces us to have $ \abs{\name EX {R''} }\ge \tfrac \kappa 4  {\theta} $. 
And this implies that $ \rho \in \mathcal S _{R''}$, as in (\ref{e.uni1}).  This is 
the desired contradiction.

\end{proof}

Our third and final fact about the collection $ \mathcal U(R)$ is a consequence 
of Lemma~\ref{l.samedirection} and a geometric observation of 
J.-O.~Stromberg \cite[Lemma 2, p. 400]{MR0481883}.

\begin{lemma}\label{l.stromberg}  For any interval $ I\subset I _{R}$ we have the inequality 
\begin{equation}\label{e.removeafew}
\sum _{\substack{R'\in \mathcal U(R)\\ \name L {R'}\le \abs I \le 
\sqrt\kappa \name L {R'} }}
\abs{ R'\cap I\times J } \le 5\abs I \cdot \name W \rho  \,.
\end{equation}
\end{lemma}

\def\thinrect#1#2#3#4{
	\pgfrect[stroke]{\pgfpoint{#3cm}{#4cm}}{\pgfpoint{#2cm}{.3cm}}
	 }

 \begin{figure}
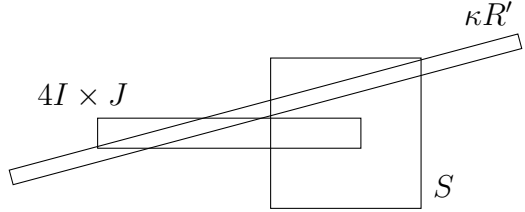
   
\begin{center} 
\begin{pgfpicture}{0cm}{0cm}{5cm}{4cm}
\begin{pgftranslate}{\pgfxy(1.5,0)}
\pgfputat{\pgfxy(-1.5,1.2)}{\pgfbox[center,center]{$4I\times J$}}
\pgfrect[stroke]{\pgfxy(-1.3,0.5)}{\pgfxy(3.5,.4)}
\pgfrect[stroke]{\pgfxy(1.0,-0.3)}{\pgfxy(2,2)}
\pgfputat{\pgfxy(3.3,0)}{\pgfbox[center,center]{$S$}}
\begin{pgftranslate}{\pgfxy(-1,.5)}
\begin{pgfrotateby}{\pgfdegree{15}}
\pgfrect[stroke]{\pgfxy(-1.5,-.1)}{\pgfxy(7,.2)}
\end{pgfrotateby}
\end{pgftranslate}
\pgfputat{\pgfxy(3.9,2.3)}{\pgfbox[center,center]{$\kappa R'$}}
\end{pgftranslate}
\end{pgfpicture}
\end{center}
\caption{The proof of Lemma~\ref{l.stromberg}} 
\end{figure}

\begin{proof}  
For each point $ x\in 4 I\times J$, consider the square $ S$ centered at $ x$ of 
side length equal to $  \sqrt \kappa \cdot  \abs I 
\cdot \abs {\operatorname e _R - \operatorname e _\rho }$.
See Figure~6. 
It is Stromberg's observation 
that for $ R'\in \mathcal U(R)$ we have 
\begin{equation*}
\frac{ \abs{\kappa R'\cap I\times J } }{\abs{ I\times J}} \simeq  \frac {\abs{S\cap \kappa R' }}{\abs {S} }
\end{equation*}
with the implied constant being independent of $ \kappa $. 
Indeed, by Lemma~\ref{l.samedirection}, we have that 
\begin{align*}
\frac{ \abs{\kappa R'\cap I\times J } }{\abs{ I\times J}}
& \simeq  \frac{ \kappa \name W {R'} }
	{\abs{ \operatorname e _{R}-\operatorname e _{\rho } }\cdot \abs I} 
\\
& \simeq \frac{ \kappa \name W {R'} \cdot \abs I 
\cdot \abs {\operatorname e _R - \operatorname e _\rho }}
{(\abs{ \operatorname e _{R}-\operatorname e _{\rho } }\cdot \abs I)^2}  
\\
&\simeq \frac {\abs{S\cap \kappa R' }}{\abs {S} }\,,
\end{align*}
as claimed.

Now, assume that (\ref{e.removeafew}) does not hold and seek a contradiction.  
Let $ \mathcal U'\subset \mathcal U(R)$ denote the collection of rectangles $ R'$ 
over which the sum is made in (\ref{e.removeafew}).  The rectangles in $ \mathcal U'$ 
were added in some order to the collection $ \mathcal R'$, and in particular 
there is a rectangle $ R_0\in \mathcal U'$ that was the last to be added to 
$ \mathcal U'$.  Let $ \mathcal U''$ be the 
collection $ \mathcal U'-\{R_0\}$.  By construction, $ \mathcal U'$ must consist of at 
least three rectangles, so that $ \mathcal U''$ is not empty.  Moreover, we certainly have 
\begin{equation*}
\sum _{R'\in \mathcal U''} \abs{ R''\cap I\times J}\ge{}4\abs{ I\times J}. 
\end{equation*}

Since we cannot have $ \rho \in \mathcal S _{R_0}$, Stromberg's observation implies that 
\begin{equation*}
R_0\subset\Bigl\{ \operatorname M _{\kappa } \sum _{R'\in \mathcal U''} \mathbf 1 _{R'}> 
\kappa ^{-
1}\Bigr\}\,.
\end{equation*}
Here, we rely upon the fact that the maximal function  $ \operatorname M _{\kappa }$ is larger than the usual maximal 
function over squares. 
But this is a contradiction to our construction, and so the proof is complete.  
\end{proof}

The principal line of reasoning to prove (\ref{e.uni3}) can now begin. 
We  make a recursive decomposition of the collection $ \mathcal U(R)$, which 
is indexed by a collection of intervals  $ \mathcal I $ that we now define.  Initialize 
\begin{equation}\label{e.principal}
\begin{split}
\mathsf{STOCK}&\leftarrow \mathcal U(R)
\\
\mathcal I&\leftarrow \emptyset 
\end{split}
\end{equation}
While there is an interval $ I\subset I_R $ satisfying 
\begin{equation*}
\sum _{\substack{R'\in \mathsf{STOCK}\\ \name L R \ge 8\abs I }} \abs{ R'\cap I \times J_R } 
\ge 10\abs I \cdot \name W \rho  \,, 
\end{equation*}
we take $ I$ to be a maximal interval with this property, and update 
\begin{gather*}
\mathcal I\leftarrow \mathcal I\cup \{I\}\,;
\\
\mathcal V (I)\leftarrow \{ R'\in \mathsf{STOCK}\mid \name L {R'} \ge 8\abs I\,, 
R'\cap I \times J_R\neq\emptyset \} \,;
\\
\mathsf{STOCK}\leftarrow \mathsf{STOCK}-\mathcal V (I)\,;
\end{gather*}
When the while loop terminates, we set $ \mathcal V\leftarrow \mathsf{STOCK}$. 

\smallskip 

It is then clear that we must have 
\begin{equation*}
\sum _{R'\in \mathcal V} \abs{ R'\cap \rho } 
\lesssim  \abs {I_R} \cdot \name W \rho \le \mathsf L (R) \cdot \mathsf W (\rho ) \,.
\end{equation*}
Lemma~\ref{l.stromberg}  implies that each $ I \in \mathcal I$ must 
have length $ \abs{ I}\le \kappa ^{-1/2} \abs{ I_\rho }$. 
But we choose intervals in $ I$ to be of maximal length, so we have 
\begin{equation} \label{e.V}
\sum _{R'\in \mathcal V(I)} \abs{ R'\cap \rho } 
\le 20\cdot  \abs I \cdot \name W \rho  \,.
\end{equation}

\medskip 
 Lemmas~\ref{l.geo}  and~\ref{l.stromberg} place significant restrictions on  
the collection of intervals $ \mathcal I$.  If we have $ I\neq I'\in \mathcal I$
with $ 2I\cap 2 I'\neq\emptyset$, then we must have e.g.~$\sqrt \kappa \abs{ I'}<\abs{ I} $.  
Moreover, we cannot have three distinct intervals $ I,I',I''\in \mathcal I$ with non 
empty intersection.  Therefore, we must have 
\begin{equation*}
\sum _{I\in \mathcal I} \abs{ I} \lesssim \abs{ I_R} \lesssim \mathsf L (R). 
\end{equation*}
With (\ref{e.V}), this completes the proof of (\ref{e.uni3}).


\end{document}